\documentclass[12pt,oneside]{amsproc}
\usepackage[T2A]{fontenc}
\usepackage[cp1251]{inputenc}
\usepackage[russian]{babel}
\usepackage{amssymb}
\usepackage{amsthm}

\title{О произведении ядерных операторов}
\author{О. И. Рейнов}
\address{Санкт-Петербургский Государственный Университет
}
\email{orein51@mail.ru}

\begin{document}


  \maketitle



 Эта работа появилась благодаря следующему вопросу Б. С. Митягина, 
 заданному в 2014 г.
на конференции, посвященной памяти А. Пелчинского, в Бедлево (Польша):
верно ли, что произведение двух ядерных операторов в банаховых пространствах 
факторизуется через 
ядерный оператор в гильбертовом пространстве?
Отметим, что этот вопрос был мотивирован анализом спектров 
компактных операторов,
некоторые степени которых ядерны, проведенным Б. С. Митягиным (см. [2]).
Мы получили сначала отрицательный ответ 
(см. ниже доказательство части следствия 2),    
а затем рассмотрели более общую ситуацию. Именно, пусть оператор $T$ 
представляет собой
композицию $T_1T_2\cdots T_n$ ядерных операторов таких, что для 
$j=1, 2, \dots, n$ оператор $T_j$ является $s_j$-ядерным, где  $0<s_j\le 1.$
Каков {\it точный}\ показатель $r=r(s_1, s_2, \dots, s_n)$, для которого $T$
факторизуется через $S_r$- оператор в гильбертовом пространстве? Здесь 
$S_r$ --- класс фон Неймана-Шаттена. Ниже мы приводим ответы на этот вопрос.

Везде далее через $X, Y, \dots$ обозначаются банаховы пространства, 
$H$ --- гильбертово пространство, $L(X, Y)$ --- банахово пространство всех
линейных непрерывных операторов из $X$  в $Y.$ 

Напомним, что оператор $T:X\to Y$ называется  {\it $s$-ядерным}\ $(0<s\le1,$ 
см., например, [4]),
      если он представим в виде
   $ Tx=\sum_{k=1}^\infty \langle x'_k,x\rangle y_k$
для 
$x\in X,$ где  $(x'_k)\subset X^*, (y_k)\subset Y,\, 
\sum_k ||x'_k||^s\,||y_k||^s<\infty.$ 
Мы используем обозначение
 $N_s(X,Y)$ для линейного пространства всех таких операторов
 и $\nu_s(T)$ для соответствующей квазинормы
$\inf   (\sum_k ||x'_k||^s\,||y_k||^s)^{1/s}.$
 В случае, когда $s=1,$
 эти операторы называют просто {\it ядерными}.

Класс  $S_p, 0<p<\infty,$ операторов фон Неймана-Шаттена определяется 
следующим образом. Пусть
$U$ --- компактный оператор в гильбертовом пространстве 
$H$ и $(\mu_n)$ --- 
последовательность его сингулярных чисел (см. [7], Теорема VI.17, с. 227).
Оператор  $U$ принадлежит пространству $S_p(H),$
если сходится ряд $\sum \mu_n^p$ (см., например, [6], 15.5.1). Пространство 
$S_p(H)$ имеет естественную квазинорму
$\sigma_p(U)=(\sum \mu_n^p)^{1/p}.$
 Отметим, что для $s\in (0,1]$ имеет место равенство
$N_s(H)=S_s(H)$ [3].        
Известно, что последовательность всех собственных чисел
оператора $U\in S_p(H),$
взятых с учетом кратностей, лежит в пространстве $l_p$ [6, теорема 27.4.3].    
Мы также будем использовать включение $S_p\circ S_q\subset S_r,$
где $0< p,q,r<\infty, 1/r=1/p+1/q$ [6, теорема 15.5.9].     

Определение идеала $\Pi_2$
абсолютно 2-суммирующих операторов можно найти в [6]  
(оно здесь нам не нужно). Отметим лишь, что в гильбертовом случае
$\Pi_2(H)= S_2(H)$ (операторы Гильберта-Шмидта; см. [6],  теорема 17.5.3).  
Кроме того, $\Pi_2(C(K), H)= L(C(K), H)$
 (см., например, [6]).    
    \smallskip
    
    {\bf Определение}.\,
Оператор
 $T$ {\it факторизуется
через оператор}\,  из $S_p(H),$ если существуют такие операторы
$A\in L(X, H),\, U\in S_p(H)$ и $B\in L(H, Y),$ что $T=BUA.$  
Если  $T$  факторизуется
через оператор  из $S_p(H),$ то  полагаем
$\gamma_{S_p}(T)= \inf ||A||\, \sigma_p(U)\, ||B||,$
где инфимум берется по всем возможным факторизациям оператора $T$
через оператор  из $S_p(H).$
    
    Имеет место следующая
    
    {\bf  Теорема}.\
    {\it
Если $X_1, X_2, \dots, X_{n+1}$ --- банаховы пространства, 
$s_k\in (0,1]$ и $T_k\in N_{s_k}(X_k, X_{k+1})$ для
$k=1, 2, \dots, n,$ то произведение   
$T:=T_n T_{n-1}\cdots T_1$
факторизуется через оператор  из $S_r(H),$
где
$1/r= 1/s_1 + 1/s_2 + \dots + 1/s_n - (n+1)/2.$      
При этом, 
$\gamma_{S_r}(T)\le \prod_{k=1}^n \nu_{s_k}(T_k).$
Результат точен.
}
\smallskip

 {\bf  Следствие 1}.\  
 {\it
Пусть $s, q\in (0,1].$
Произведение  $s$-ядерного и  $q$-ядерного операторов 
в банаховых пространствах
факторизуется через   $S_r(H)$-оператор, где
$1/r=1/s+1/q -3/2.$
Результат точен.
}
\smallskip

 {\bf  Следствие 2}.\
 {\it
Произведение  двух ядерных операторов в банаховых пространствах
факторизуется через   $S_2(H)$-оператор. Результат точен.
}
\smallskip

Проведем доказательство (эскиз) второй части следствия 2,
из которого получится и отрицательный ответ на упомянутый выше 
вопрос Б. С. Митягина.

Пусть $f$ --- непрерывная на единичной окружности $\mathbf T$
функция   Карлемана [1] (см. также [5], теорема 4.11, с. 321),   
последовательность коэффициентов Фурье 
которой лежит в 
$ l_2\setminus \cup_{p<2} l_p.$ Рассмотрим на $C:= C(\mathbf T)$ 
оператор свертки с функцией $f:$
$$T: C\overset{*f}\to C.$$
Оператор $T$ является ядерным (как интегральный оператор с непрерывным ядром).
Рассмотрим произведение $TT.$ Отметим, что последовательность 
собственных чисел этого 
оператора лежит в $l_1,$ но не в $l_s, 0<s<1.$
Пусть $r\in [1, 2].$ Предположим, что существует такой 
$S_r$-оператор $U\in S_r(H),$ что $TT$ факторизуется следующим образом:
$$TT: C\overset A\to H\overset U \to H\overset B\to C.$$
Рассмотрим диаграмму
$$H\overset B\to C\overset A\to H\overset U \to H\overset B\to C.$$
Набор собственных чисел оператора 
$UAB$ совпадает (с учетом их алгебраических кратностей) 
с набором всех собственных
чисел оператора  $TT=BUA$ [6, стр. 436]    
(и, таким образом, лежит в $l_1$ и не лучше).  Но:
$$A\in \Pi_2;\  \text{следовательно, }\ AB\in S_2;\  U\in S_r.$$
Отсюда: $UAB\in S_{s},$
где $1/s=1/2+ 1/r,$
и, таким образом, последовательность собственных чисел этого оператора лежит в
$l_s.$  Следовательно, $s=1$  и $r=2.$  
\smallskip


 \newpage
 
 \centerline{ЛИТЕРАТУРА}
  
   \bigskip
   
            [{\bf 1}]\,
   T. Carleman,
   Acta Mathematica, {\bf 41}:1 (1916), 377-384.\
               [{\bf 2}]\,     
   B. S. Mityagin,
{J. Operator Theory}, \textbf{76}:1 (2016), 57--65.\
             [{\bf 3}]\,
   R. Oloff,  
   Beitr\"age Anal., {\bf 4}  (1972), 105-108.\
   [{\bf 4}]\,     
    O. I. Reinov,
 J.  Math. Sciences, \textbf{115}:3 (2003), 2243-2250.\
           [{\bf 5}]\,   
  А. Зигмунд,     
  {\it Тригонометрические ряды}\ I,\,
  Мир, Москва, 1965.  
            [{\bf 6}]\,   
   А. Пич,
     {\it Операторные идеалы},
  Мир, Москва, 1982.   
             [{\bf 7}]\,   
М. Рид,  Б. Саймон,
 {\it Методы современной математической физики}\ 1,
 {\it Функциональный анализ}, 
   Мир, Москва, 1977. 
       
\bigskip

{\bf  Замечание} (добавлено после появления заметки в печати).\  
Точность первого утверждения теоремы может быть получена методом,
аналогичным доказательству следствия 2: вместо примера Карлемана надо
рассматривать аналогичные примеры функций из подходящих классов Липшица
(см., например, [5]). 
Соответствующие результаты, возможно, будут опубликованы.

В нашей статье “О произведении s-ядерных операторов” 
(Матем. заметки, 107:2 (2020), 311-316)
мы пошли иным путем: привели конечномерные аналоги теоремы
 и применили их к доказательству этой теоремы в полной общности.      
       
\end{document}